\documentclass{article}
\usepackage{amsmath}
\usepackage{amssymb}
\usepackage[amsmath, thmmarks]{ntheorem}

\newcommand{\C}{\mathbb{C}}
\newcommand{\PP}{\mathbb{P}}

\newcommand{\Q}{\mathbb{Q}}
\newcommand{\pp}{\mathfrak{p}}
\newcommand{\FF}{\mathcal{F}}

\newcommand{\ord}{\mbox{ord}}
\newcommand{\Div}{\mbox{div}}
\newcommand{\D}{\mathbb{D}}

\theoremseparator{.}

\theoremstyle{change}

\theoremstyle{nonumberplain}
\newtheorem{cor*}{Corollary}
\renewtheorem{lem*}{Lemma}

\renewtheorem{prop*}{Proposition}

\theorembodyfont{\upshape} 
\theoremstyle{nonumberplain}
\theoremheaderfont{\itshape}
\theoremsymbol{\ensuremath{\Box}}
\newtheorem{proof}{Proof} \theoremsymbol{}
\newtheorem{rem}{Remark}

\begin{document}
\title{A remark on a relation between foliations and number theory}
\author{Fabian Kopei}
\date{}
\maketitle
\begin{abstract}
We interpret a formula for meromorphic functions on
foliations by Riemann surfaces as an analogue to the product
formula of valuations in algebraic number theory. 
\end{abstract}
A meromorphic function on a compact Riemann surface has as
many zeros as poles. This well-known geometric fact has an algebraic analogue:
for a smooth algebraic curve $C$ and an element $f$ of its function
field we have
\begin{eqnarray}\sum_{P\in C} \ord_P f=\deg(\Div(f))=0.\label{picard}\end{eqnarray}

In the algebraic context a similar formula holds for the
valuations of an algebraic number field (i.e.\! a finite extensions of $\Q$) 
or a separable function field:

$$\sum_{\pp \mbox{ finite}}\log |f|_\pp\
+\negthickspace\negthickspace\sum_{\pp \mbox{ infinite}}
\negthickspace\negthickspace\log |f|_\pp= -\negthickspace\negthickspace\sum_{\pp \mbox{ prime}}\ord_\pp\,(f)\log\mathfrak{N}(\pp)\,+\negthickspace\negthickspace\sum_{\pp \mbox{ infinite}}
\negthickspace\negthickspace\log |f|_\pp=0,$$
where $f$ is a non-zero element of the field and the sums run over the equivalence classes of the finite or
infinite valuations, respectively the prime-ideals of the ring of integers.

In this article we want to search a geometric
analogue for this formula.
Neglecting for the moment the infinite primes, the structural
difference to equation \nolinebreak (\ref{picard}) is the
factor $\log\mathfrak{N}(\pp)$. In an expected geometric analogon this
is only a property of the point $P$, not of the function $f$. Since it 
seems to be difficult to assign real values to points in a canonical
way, we may try to replace the point by the next simple
geometric object, a line, to which we assign its length measured by a
flow. But a meromorphic function has only isolated zeros and poles, so
this line should come with a transverse structure. The simplest
object realizing this seems to be a 3-dimensional foliation by
Riemann surfaces. 
Assuming that
the flow and the foliation are compatible with each other, we may
therefore expect the following

\begin{prop*}
Let $M$ be an oriented, closed 3-manifold, $\FF$ a leafwise oriented
foliation by Riemann surfaces and $\Phi$ a transverse foliation-invariant
flow. Let  
$$f:M\longrightarrow \C\PP^1$$ be a smooth function which is
meromorphic on the leaves, such that the zeros and poles lie only on
finitely many closed orbits $\gamma_1,\dots,\gamma_n$.

Then on these
orbits $\ord_{\gamma_i}f$ is constant and we have
$$\sum_{i=1}^n l(\gamma_i)\ord_{\gamma_i}f=0,$$
where $l(\gamma_i)$ denotes the length of the orbit, i.e. the smallest 
return time.
\end{prop*}

\begin{rem}
The proposition is a special case of a well-known theorem for
laminations: the harmonic measure of zeros minus that of poles is equal to
the integrated Chern class, see \cite{Gh}. For similar formulas see \cite{Co}, \cite{Hu}.  
\end{rem}

\begin{rem}
For a more general definition of meromorphic functions see for example 
\cite{Gh}, concerning their existence see \cite{Gh}, \cite{De}.
\end{rem}

So for the case of no infinite primes we found a geometric analogue:
here elements of the field correspond to meromorphic functions modulo $\C^*$ and 
primes correspond to closed orbits with length $\log\mathfrak{N}(\pp)$.

The infinite primes are something like a compactification of the finite
ones. For an analogue we may therefore try to take a noncompact
space, which we then compactify. The simplest example is, that we have
to complete our foliation by a finite number of compact leaves, on
which the flow is not transverse. We get then an analogue for the case of
infinite primes:
\begin{prop*}
Let $M$ be an oriented, closed 3-manifold, $\FF$ a leafwise oriented
foliation by Riemann surfaces and $\Phi$ an foliation-invariant
flow. We assume that the flow is up to a finite number of compact
leaves $\left(L_j\right)_{j=1,\dots,k}$ transverse. Let $$f:M\setminus \cup_j L_j\to \C\PP^1$$ be a smooth function
which is  
meromorphic on the leaves, such that the zeros and poles lie only on
finitely many closed, transverse orbits $\gamma_1,\dots,\gamma_n$. Let $\omega$ be the on
$M\setminus \cup_j L_j$ defined
1-form, which is zero on the integrable distribution of the foliation and one on the vector field
generated by the flow. We set
$$\eta:=\frac{1}{f}d_\FF f\wedge\omega,$$
where $d_\FF$ is the differential along the leaves. Let $TL_j$ be
pairwise distinct compact tubular neighbourhoods of the leaves $L_j$, such
that the $\gamma_i$ and the $TL_j$ are pairwise distinct. Then
$$\sum_{i=1}^n l(\gamma_i)\ord_{\gamma_i}f=-\sum_{j=1}^k\frac{1}{2\pi
i}\int_{\partial TL_j}\eta.$$
\end{prop*}
The proof is a simple adaption of \cite{Gh}.
\begin{proof}It is an easy application of the theorem of Cauchy, that
$\ord_{\gamma_i}f$ is constant on the orbits.

For each $\gamma_i$ we take distinct, orientation preserving tubular neighbourhoods
$$\iota_i:S^1\times \D\hookrightarrow M\setminus \cup_j L_j,$$
where $\D\subset \C$ is the 2-dimensional disc, which respect the foliation, i.e. $\iota(\{t\}\times \D)$
lies in a leaf (\cite{CC}, p. 89). We have
$$\int_{\iota_i(S^1\times \partial \D)}\eta=\int_{S^1\times \partial \D}(\iota_i^{-1})^*(\eta)=\int_0^{l(\gamma_i)}dt
\oint\frac{f'(\iota_i(z,t))}{f(\iota_i(z,t))}dz=2\pi i\ l(\gamma_i)\ord_{\gamma_i} f,$$
where the last equation is a consequence of the theorem of the zeros and poles
counting integral. 

So if $\Gamma_i$ denotes the image of $\iota_i$, we have 
on the one hand 
$$\sum_{i=1}^n \int_{\partial \Gamma_i}\eta=2\pi i\sum_{i=1}^n l(\gamma_i)\ord_{\gamma_i}f.$$
On the other hand by the theorem of Stokes and the fact that $\eta$
is closed we have
$$\sum_{i=1}^n \int_{\partial \Gamma_i}\eta+\sum_{j=1}^k
\int_{\partial TL_j}\eta=-\int_{\partial \left(M\setminus \left( \cup_i
\Gamma_i \cup_j TL_j\right)\right)}\eta = -\int_{M\setminus\left( \cup_i
\Gamma_i\cup_j TL_j\right)}d\eta= 0.$$
The combination of these two equations proves the proposition.
\end{proof}

\begin{rem}In this article we wanted to point out a similarity between
certain foliations and number theory in a self-contained way. But of
course this work is strongly influenced by an elaborate
conjecture, in which Deninger postulates a deep relation between these 
objects (see for example \cite{De01}). The work of Deninger suggests
also a correspondence between infinite primes and compact,
non-transversal leaves.
\end{rem}

Mathematisches Institut\\
Westf. Wilhelms-Universit\"at\\
Einsteinstr. 62\\
48149 M\"unster\\
Germany\\
f\_kope01@math.uni-muenster.de
   
\end{document}